\documentclass[]{article}
\usepackage{proceed2e}

\usepackage{times}
\usepackage{pifont}
\usepackage[T1]{fontenc}
\usepackage{currvita}

\usepackage{pdfpages}
\usepackage{tikz}
\usetikzlibrary{arrows}
\usetikzlibrary{positioning}

\usepackage{amsthm}
\usepackage{amsmath}
\usepackage{amssymb}
\usepackage[authoryear]{natbib}

\usepackage{caption}
\usepackage{subcaption}
\usepackage{graphicx}

\usepackage{xr}
\RequirePackage{hyperref}

\DeclareMathOperator{\CPDAG}{CPDAG}
\DeclareMathOperator{\DAG}{DAG}
\DeclareMathOperator{\MAG}{MAG}
\DeclareMathOperator{\PAG}{PAG}
\DeclareMathOperator{\PossDe}{PossDe}
\DeclareMathOperator{\PossPossDe}{(Poss)De}

\DeclareMathOperator{\De}{De}
\DeclareMathOperator{\Pa}{Pa}
\DeclareMathOperator{\distancefrom}{distance-from-\!}

\DeclareMathAlphabet{\mathpzc}{OT1}{pzc}{m}{it}

\newenvironment{proofsketchamenab}{\noindent{\bf Proof of Lemma \ref{lemmaadjamengen}:}}%
{\hfill $\square$}
{\hfill $\square$}
{\hfill $\square$}
{\hfill $\square$}
{\hfill $\square$}

\title{A Complete Generalized Adjustment Criterion}

\author{ {\bf Emilija Perkovi\'c} \\
Seminar for Statistics \\
ETH Zurich, Switzerland\\
perkovic@stat.math.ethz.ch \\
\And
{\bf Johannes Textor}\\
Theoretical Biology \& Bioinformatics \\
Utrecht University, The Netherlands\\
johannes.textor@gmx.de \\
\And
{\bf Markus Kalisch}\\
Seminar for Statistics \\
ETH Zurich, Switzerland\\
kalisch@stat.math.ethz.ch\\
\And
{\bf Marloes H. Maathuis}\\
Seminar for Statistics \\
ETH Zurich, Switzerland\\
maathuis@stat.math.ethz.ch\\
}


\renewcommand{\bibname}{References}
\newtheorem{theorem}{Theorem}[section]

\newtheorem{lemma}[theorem]{Lemma}
\newtheorem{definition}[theorem]{Definition}
\newtheorem{example}[theorem]{Example}

\theoremstyle{definition}

{\hfill $\square$}
{\hfill $\square$}
\newenvironment{proofthm}{\noindent{\bf Proof of Theorem \ref{thmgenac}:}}%
{\hfill $\square$}
\newenvironment{prooflemmaeqa}{\noindent{\bf Proof of Lemma \ref{lemmaprovingeqofacconditiona}:}}%
{\hfill $\square$}
\newenvironment{prooflemmaeqb}{\noindent{\bf Proof of Lemma \ref{lemmaequivalenceofcondb}:}}%
{\hfill $\square$}
{\hfill $\square$}
\newenvironment{proofof}[1][]{\begin{trivlist}
   \item[\hskip \labelsep {\bfseries Proof of #1.}]}{\end{trivlist} \hfill $\square$}

\theoremstyle{remark}

\newcommand{\vars}[1][V]{\mathbf{#1}}
\newcommand{\e}[1][E]{\mathbf{#1}}

 \newcommand{\g}[1][G]{\mathcal{#1}}

\newcommand{\f}[2][X,Y]{\mathbf{F}_{\mathpzc{#2}}(#1)}
\newcommand{\fb}[2][X,Y]{\mathbf{F}_{#2}(\mathbf{#1})}

\newcommand{\pstar}[1][p]{{#1}^{*}}

\newcommand{\bulletcirc}{
  \setlength{\unitlength}{1mm}
  \begin{picture}(5,1)(0,0)
    \put(0.2,0){$\bullet$}
    \put(1.1, 1){\line(1,0){2.4}}
    \put(4, 1){\circle{1}}
  \end{picture}
}

\newcommand{\circbullet}{
  \setlength{\unitlength}{1mm}
  \begin{picture}(5,1)(0,0)
    \put(1,1){\circle{1}}
    \put(1.5,1){\line(1,0){2.4}}
    \put(2.9,0){$\bullet$}
  \end{picture}
}

\newcommand{\bulletarrow}{
  \setlength{\unitlength}{1mm}
  \begin{picture}(5,1)(0,0)
    \put(0.2,0){$\bullet$}
    \put(1,0){$\rightarrow$}
  \end{picture}
}

\newcommand{\arrowbullet}{
  \setlength{\unitlength}{1mm}
  \begin{picture}(5,1)(0,0)
    \put(0.2,0){$\leftarrow$}
    \put(3,0){$\bullet$}
  \end{picture}
}

\newcommand{\circarrow}{
  \setlength{\unitlength}{1mm}
  \begin{picture}(5,1)(0,0)
    \put(1,1){\circle{1}}
    \put(1.2,0){$\rightarrow$}
  \end{picture}
}

\newcommand{\arrowcirc}{
  \setlength{\unitlength}{1mm}
  \begin{picture}(5,1)(0,0)
    \put(0.2,0){$\leftarrow$}
    \put(4.3,1){\circle{1}}
  \end{picture}
}

\newcommand{\circcirc}{
  \setlength{\unitlength}{1mm}
  \begin{picture}(5,1)(0,0)
    \put(1,1){\circle{1}}
    \put(1.5,1){\line(1,0){2}}
    \put(4,1){\circle{1}}
  \end{picture}
}

\begin{document}

\maketitle

\begin{abstract}

Covariate adjustment is a widely used
approach to estimate total causal effects from observational
data. Several graphical criteria have been developed
in recent years to identify valid covariates for adjustment
from graphical causal models.
These criteria can handle multiple
causes, latent confounding, or partial knowledge
of the causal structure;
however, their diversity is confusing
and some of them are only sufficient, but not necessary.
In this paper, we present a criterion that is
necessary and sufficient
for four different classes of graphical causal
models: directed acyclic graphs ($\DAG$s), maximum ancestral
graphs  ($\MAG$s), completed partially directed acyclic graphs
($\CPDAG$s), and partial ancestral graphs ($\PAG$s).
Our criterion subsumes the existing ones and in this way
unifies adjustment set construction for a large set of graph classes.
\end{abstract}

\section{INTRODUCTION}\label{sec: intro}

Which covariates do we need to adjust for when estimating
total causal effects from observational data? Graphical
causal modeling allows to answer this question
constructively, and contributed fundamental insights
to the theory of adjustment in general.
For instance, a simple example known as the ``M-bias graph''
shows that it is not always appropriate to
adjust for all observed (pre-treatment) covariates
\citep{Shrier2008,Rubin2008}.
A few small graphs also suffice to refute the ``Table 2 fallacy''
\citep{Westreich2013}, which is the belief that the
coefficients in multiple regression models are ``mutually
adjusted''. Thus,
causal graphs had substantial impact on
theory and practice of covariate adjustment \citep{ShrierP2008}.

The practical importance of covariate adjustment has inspired
a growing body of theoretical work on graphical criteria that are sufficient and/or necessary for adjustment.
Pearl's back-door criterion \citep{pearl1993bayesian} is probably the most well-known, and is
sufficient but not necessary for adjustment in $\DAG$s.
\citet{shpitser2012validity} adapted the back-door criterion to a necessary and sufficient graphical criterion
for adjustment in $\DAG$s.
Others considered graph classes other than
$\DAG$s, which can represent structural uncertainty.
\citet{vanconstructing} gave necessary and sufficient graphical criteria for $\MAG$s that allow for unobserved variables (latent confounding). \citet{maathuis2013generalized}
presented a generalized back-door criterion for $\DAG$s, $\MAG$s, $\CPDAG$s and $\PAG$s,
where $\CPDAG$s and $\PAG$s represent Markov equivalence classes of $\DAG$s or $\MAG$s, respectively,
and can be inferred directly from data (see, e.g., \citealp{spirtes2000causation, Chickering03, Colombo2012, ClaassenMooijHeskes13, colombo2014order}). The generalized back-door criterion is sufficient but not necessary for adjustment.

In this paper, we extend the results of
\citet{shpitser2012validity}, \citet{vanconstructing} and \citet{maathuis2013generalized} to derive a single necessary and sufficient adjustment criterion that holds for all four graph classes: $\DAG$s,
$\CPDAG$s, $\MAG$s and $\PAG$s.

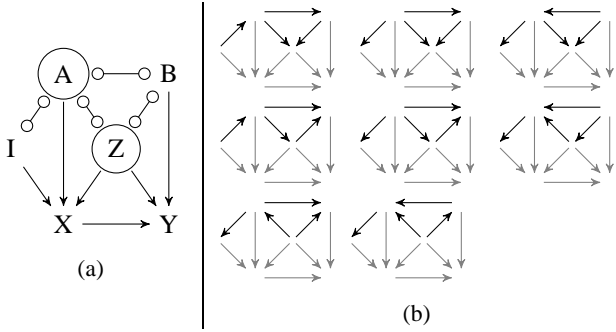
\begin{figure}

\tikzstyle{every edge}=[draw,>=stealth',->]
\newcommand\dagvariant[1]{\begin{tikzpicture}[xscale=.5,yscale=0.5]
\node (a) at (0,2) {};
\node (b) at (2,2) {};
\node [inner sep=1pt] (i) at (-1,1) {};
\node [inner sep=1pt] (z) at (1,1) {};
\node (x) at (0,0) {};
\node (y) at (2,0) {};
\begin{scope}[gray]
\draw (a) edge (x);
\draw (b) edge (y);
\draw (z) edge (x) edge (y);
\draw (x) edge (y);
\draw (i) edge (x);
\end{scope}
\draw #1;
\end{tikzpicture}}

\begin{subfigure}{.3\columnwidth}

\begin{tikzpicture}[xscale=.7,yscale=1]
\node [circle,draw,outer sep=1pt] (a) at (0,2) {A};
\node (b) at (2,2) {B};
\node (i) at (-1,1) {I};
\node [circle,draw,outer sep=1pt] (z) at (1,1) {Z};
\node (x) at (0,0) {X};
\node (y) at (2,0) {Y};
\draw (i) edge [o-o] (a) edge  (x);
\draw (z) edge [o-o] (a) edge [o-o] (b) (a) edge [o-o] (b);
\draw (a) edge (x);
\draw (b) edge (y);
\draw (z) edge (x) edge (y);
\draw (x) edge (y);
\end{tikzpicture}

\caption{}
\end{subfigure}
\hfill
\vrule
\hfill
\begin{subfigure}{.65\columnwidth}
\dagvariant{
(i) edge [->] (a)
(z) edge [<-] (a) edge [<-] (b)
(a) edge [->] (b)
}
\dagvariant{
(i) edge [<-] (a)
(z) edge [<-] (a) edge [<-] (b) (a) edge [->] (b)
}
\dagvariant{
(i) edge [<-] (a)
(z) edge [<-] (a) edge [<-] (b) (a) edge [<-] (b)
}
\dagvariant{
(i) edge [->] (a)
(z) edge [<-] (a) edge [->] (b) (a) edge [->] (b)
}
\dagvariant{
(i) edge [<-] (a)
(z) edge [<-] (a) edge [->] (b) (a) edge [->] (b)
}
\dagvariant{
(i) edge [<-] (a)
(z) edge [->] (a) edge [<-] (b) (a) edge [<-] (b)
}
\dagvariant{
(i) edge [<-] (a)
(z) edge [->] (a) edge [->] (b) (a) edge [->] (b)
}
\dagvariant{
(i) edge [<-] (a)
(z) edge [->] (a) edge [->] (b) (a) edge [<-] (b)
}
\caption{}
\end{subfigure}

\caption{(a) A $\CPDAG$ in which, according to our
criterion, $\{A,Z\}$ is an  adjustment set
for the total causal effect of $X$ on $Y$. (b)
The Markov equivalence class of (a), with node
labels removed for simplicity and varying edges
highlighted. An adjustment
set for a $\CPDAG$ ($\PAG$) is one that
works in all $\DAG$s ($\MAG$s) of the Markov equivalence class.
}
\label{fig:cpdagexample}
\end{figure}

To illustrate the use of our generalized adjustment
criterion, suppose we are given the $\CPDAG$
in Figure~\ref{fig:cpdagexample}a and we want to
estimate the total causal effect of $X$ on $Y$.
Our criterion will inform us that the set $\{A,Z\}$ is an
adjustment set for this $\CPDAG$, which means that it
is an adjustment set in every $\DAG$ that the $\CPDAG$
represents (Figure~\ref{fig:cpdagexample}b). Hence, we can estimate the causal
effect without knowledge of the full causal structure.
In a similar manner, by applying our criterion to a $\MAG$
or a $\PAG$, we find adjustment sets that are
valid for all $\DAG$s represented by this $\MAG$ or $\PAG$.
Our criterion finds such
adjustment sets whenever they exist;
else, our knowledge of
the model structure is insufficient to compute the
desired causal effect by covariate adjustment.
We hope that this ability to allow for incomplete
structural knowledge or latent confounding or
both will help address concerns
that graphical causal modelling ``assumes that all [...] $\DAG$s
have been properly specified'' \citep{West2014}.

We note that, although we can find all causal effects that are identifiable by covariate adjustment, we generally do not find
all identifiable causal effects, since some effects may be identifiable by other means, such as
Pearl's front-door criterion \citep[Section 3.3.2]{Pearl2009} or the ID algorithm \citep{tian2002general,shpitser2006identification}.
We also point out that $\MAG$s and $\PAG$s are in principle
not only able to represent unobserved confounding, but
can also account for unobserved selection variables.
However, in this paper we assume that there are no
unobserved selection variables. This restriction is mainly due to the fact that
selection bias often renders it impossible to
identify causal effects using
just covariate adjustment. \citet{Bareinboim2014} discuss
these problems and present creative approaches to work
around them, e.g., by combining data from different sources.
We leave the question whether adjustment
could be combined with such auxiliary methods
aside for future research.

\section{PRELIMINARIES}\label{sec: prelim}

Throughout the paper we denote sets in bold uppercase letters (e.g., $\mathbf{S}$), graphs in calligraphic font (e.g., $\g$) and nodes in a graph in uppercase letters (e.g., $X$)

\textbf{Nodes and edges.} A graph $\g= (\vars,\e) $ consists of a set of nodes (variables) $ \vars=\left\lbrace X_{1},\dots,X_{p}\right\rbrace$ and a set of edges $ \e $.

There is at most one edge between any pair of nodes, and nodes are called \textit{adjacent} if they are connected by an edge. Every edge has two edge marks that can be arrowheads, tails or circles. Edges can be \emph{directed} $\rightarrow$,  \emph{bidirected}  $\leftrightarrow$, \emph{non-directed} $\circcirc$ or \emph{partially directed} $\circarrow$. We use $\bullet$ as a stand in for any of the allowed edge marks. An edge is \textit{into} (\textit{out of}) a node $X$ if the edge has an arrowhead (tail) at $X$. A \emph{directed graph} contains only directed edges. A \emph{mixed graph} may contain directed and bi-directed edges. A \emph{partial mixed graph} may contain any of the described edges. Unless stated otherwise, all definitions apply for partial mixed graphs.

\textbf{Paths.} A \textit{path} $p$ from $X$ to $Y$ in $\g$ is a sequence of distinct nodes $\langle X, \dots,Y \rangle$ in which every pair of successive nodes is adjacent in $\g$. A node $V$ \emph{lies on a path} $p$ if $V$ occurs in the sequence of nodes.
The \textit{length} of a path equals the number of edges on the path.  A \textit{directed path} from $X$ to $Y$ is a path from $X$ to $Y$ in which all edges are directed towards $Y$, i.e., $X \to \dots\to Y$. A directed path is also called a \textit{causal path}. A \textit{possibly directed path} (\textit{possibly causal path}) from $X$ to $Y$ is a path from $X$ to $Y$ that has no arrowhead pointing to $X$. A path from $X$ to $Y$ that is not possibly causal is called a \textit{non-causal path} from $X$ to $Y$. A directed path from $X$ to $Y$ together with an edge $Y\to X$ ($Y \leftrightarrow X$) forms an \emph{(almost) directed cycle}.
For two disjoint subsets $\mathbf{X}$ and $\mathbf{Y}$ of $\mathbf{V}$, a path from $\mathbf{X}$ to $\mathbf{Y}$
is a path from some $X \in \mathbf{X}$ to some $Y \in \mathbf{Y}$.
A path from $\mathbf{X}$ to $\mathbf{Y}$ is \textit{proper} if only its first node is in $\mathbf{X}$.

\textbf{Subsequences and subpaths.} A \textit{subsequence} of a path $p$ is a sequence of nodes obtained by deleting some nodes from $p$ without changing the order of the remaining nodes. A subsequence of a path is not necessarily a path. For a path $p = \langle X_1,X_2,\dots,X_m \rangle$, the \textit{subpath} from $X_i$ to $X_k$ ($1\le i\le k\le m)$ is the path $p(X_i,X_k) = \langle X_i,X_{i+1},\dots,X_{k}\rangle$. We denote the concatenation of paths by $\oplus$, so
that for example $p = p(X_1,X_{k}) \oplus p(X_{k},X_{m})$. We use the convention that we remove any loops that may occur due to the concatenation, so that the result is again a path.

\textbf{Ancestral relationships.} If $X\to Y$, then $X$ is a \textit{parent} of $Y$. If there is a (possibly) directed path from $X$ to $Y$, then $X$ is a \textit{(possible) ancestor} of $Y$, and $Y$ is a \textit{(possible) descendant} of $X$. Every node is also a descendant and an ancestor of itself.  The sets of parents and (possible) descendants of $X$ in $\g$ are denoted by $\Pa(X,\g)$ and $\PossPossDe(X,\g)$ respectively. For a set of nodes $\mathbf{X} \subseteq \mathbf{V}$, we have $\Pa(\mathbf{X},\g) = \cup_{X \in \mathbf{X}} \Pa(X,\g)$, with analogous definitions for $\PossPossDe(\mathbf{X},\g)$.

\textbf{Colliders and shields.} If a path $p$ contains $X_i \bulletarrow X_j \arrowbullet X_k$ as a subpath, then $X_j$ is a \textit{collider} on $p$. A \textit{collider path} is a path on which every non-endpoint node is a collider. A  path of length one is a trivial collider path. A path $\langle X_{i},X_{j},X_{k} \rangle$ is an \emph{(un)shielded triple} if $ X_{i} $ and $ X_{k}$ are (not) adjacent. A path is \textit{unshielded} if all successive triples on the path are unshielded. Otherwise the path is \emph{shielded}. A node $X_{j}$ is a \textit{definite non-collider} on a path $p$ if there is at least one edge out of $X_{j}$ on $p$, or if $X_{i} \bulletcirc X_j \circbullet X_k$ is a subpath of $p$ and $\langle X_i,X_j,X_k\rangle$ is an unshielded triple. A node is of \textit{definite status} on a path if it is a collider or a definite non-collider on the path. A path $p$ is of definite status if every non-endpoint node on $p$ is of definite status. An unshielded path is always of definite status, but a
definite status path is not always unshielded.

\textbf{m-separation and m-connection.} A definite status path \textit{p} between nodes $X$ and $Y$ is \textit{m-connecting} given a set of nodes $\mathbf{Z}$ ($X,Y \notin \mathbf{Z}$) if every definite non-collider on $p$ is not in $\mathbf{Z}$, and every collider on $p$ has a descendant in $\mathbf{Z}$. Otherwise $\mathbf{Z}$ blocks $p$.  If $\mathbf{Z}$ blocks all definite status paths between $X$ and $Y$, we say that $X$ and $Y$ are m-separated given $\mathbf{Z}$. Otherwise, $X$ and $Y$ are m-connected given $\mathbf{Z}$. For pairwise disjoint subsets $\mathbf{X}$, $\mathbf{Y}$, $\mathbf{Z}$ of $\vars$, $\mathbf{X}$ and $\mathbf{Y}$ are m-separated given $\mathbf{Z}$ if $X$ and $Y$ are m-separated by $\mathbf{Z}$ for any $X \in \mathbf{X}$ and $Y\in \mathbf{Y}$. Otherwise, $\mathbf{X}$ and $\mathbf{Y}$ are m-connected given $\mathbf{Z}$.

\textbf{Causal Bayesian networks.} A directed graph without directed cycles is a \emph{directed acyclic graph ($\DAG)$}. A Bayesian network for a set of variables $\vars=\{X_{1},\dots,X_{p}\}$ is a pair ($\g,f$), where $\g$ is a $\DAG$, and $f$ is a joint probability density for $\vars$ that factorizes according to the conditional independence relationships described via m-separation, that is $f(\vars)= \prod_{i=1}^{p}f(X_{i}|Pa(X_{i},\g))$ \citep{Pearl2009}. We call a $\DAG$ causal when every edge $X_{i} \rightarrow X_{j}$ in $\g$ represents a direct causal effect of $X_{i}$ on $X_{j}$. A Bayesian network ($\g,f$) is a \emph{causal Bayesian network} if $\g$ is a causal $\DAG$. If a causal Bayesian network is given and all variables are observed one can easily derive post-intervention densities. In particular, we consider interventions $do(\mathbf{X} =\mathbf{x})$ ($\mathbf{X}\subseteq \mathbf{V}$), which represent outside interventions that set $\mathbf{X}$ to $\mathbf{x}$ \citep{Pearl2009}:
\begin{multline}
f(\mathbf{v}|do(\mathbf{X}=\mathbf{x})) =   \\
=
\begin{cases}
\prod_{X_{i} \in \vars \setminus \mathbf{X}}f(x_{i}|Pa(x_{i},\g)), &  \text{for values of }\mathbf{V} \\ & \text{consistent with }\mathbf{x,} \\
0, & \text{otherwise.}
\end{cases}
\end{multline}
Equation (1) is known as the truncated factorization formula \citep{Pearl2009} or the g-formula \citep{robins1986new}.

\textbf{Maximal ancestral graph.}
A mixed graph $\g$ without directed cycles and almost directed cycles is called \textit{ancestral}. A \emph{maximal ancestral graph ($\MAG$)} is an ancestral graph $\g = (\vars, \e)$ where every two non-adjacent nodes $X$ and $Y$ in $\g$ can be m-separated by a set $\mathbf{Z} \subseteq \vars \setminus \{X,Y\}$. A $\DAG$ with unobserved variables can be uniquely represented by a $\MAG$ that preserves the ancestral and m-separation relationships among the observed variables \citep{richardson2002ancestral}.
The $\MAG$ of a causal $\DAG$ is a \emph{causal $\MAG$}.

\textbf{Markov equivalence.}
 Several $\DAG$s can encode the same conditional independence information via m-separation. Such $\DAG$s form a \emph{Markov equivalence class} which can be described uniquely by a \emph{completed partially directed acyclic graph ($\CPDAG$)} .
 Several $\MAG$s can also encode the same conditional independence information. Such $\MAG$s form a Markov equivalence class which can be described uniquely by a \emph{partial ancestral graph ($\PAG$)} \citep{richardson2002ancestral, ali2009markov}.
We denote all $\DAG$s ($\MAG$s) in the Markov equivalence class described by a $\CPDAG$ ($\PAG$) $\g$ by $[\g]$.

\textbf{Consistent density.} A density $f$ is \textit{consistent} with a causal $\DAG$ $\g[D]$ if the pair ($\g[D],f$) forms a causal Bayesian network. A density $f$ is consistent with a causal $\MAG$ $\g[M]$ if there exists a causal Bayesian network ($\g[D]',f'$), such that $\g[M]$ represents $\g[D]'$ and $f$ is the observed marginal of $f'$. A density $f$ is consistent with a $\CPDAG$ ($\PAG$) $\g$ if it is consistent with a $\DAG$ ($\MAG$) in $[\g]$.

\textbf{Visible and invisible edges.}
All directed edges in $\DAG$s and $\CPDAG$s are said to be visible. Given a $\MAG$ $\g[M]$ or a $\PAG$ $\g$, a directed edge $X \rightarrow Y$ is \textit{visible} if there is a node $V$ not adjacent to $Y$ such that there is an edge between $V$ and $X$ that is into $X$, or if there is a collider path from $V$ to $X$ that is into $X$ and every non-endpoint node on the path is a parent of $Y$. Otherwise, $X \rightarrow Y$ is said to be \textit{invisible} \citep{zhang2006causal, maathuis2013generalized}.

\begin{figure}[!tbp]
\centering
\begin{subfigure}{.15\textwidth}
  \centering
  \begin{tikzpicture}[>=stealth',shorten >=1pt,auto,node distance=2cm,main node/.style={minimum size=0.6cm,font=\sffamily\Large\bfseries},scale=0.6,transform shape]
\node[main node]         (X)                        {$X$};
\node[main node]         (V) [above left of= X]  	{$V$};
\node[main node]       	 (Y) [right of= X]         	{$Y$};
\draw[*->] (V) edge    (X);
\draw[->] (X) edge    (Y);
\end{tikzpicture}
  \label{visible:sub1}
\end{subfigure}%
\unskip
\vrule
\begin{subfigure}{.3\textwidth}
  \centering
  \begin{tikzpicture}[>=stealth',shorten >=1pt,auto,node distance=2cm,main node/.style={minimum size=0.6cm,font=\sffamily\Large\bfseries},scale=.58,transform shape]
\node[main node]         (X)                        {$X$};
\node[main node]         (Vn) [above left of= X]  	{$V_n$};
\node[main node]         (V3) [above left of= Vn]  {$V_3$};
\node[main node]         (V2) [above left of= V3]  	{$V_2$};
\node[main node]         (V1) [above left of= V2]  	{$V_1$};
\node[main node]       	 (Y) [right of= X]         	{$Y$};
\draw[*->] (V1) edge    (V2);
\draw[<->] (V2) edge    (V3);
\draw[dotted] (V3) edge    (Vn);
\draw[<->] (Vn) edge    (X);
\draw[->] (X) edge    (Y);
\draw[->] (Vn) ..  controls (0.5,1.5) ..  (Y);
\draw[->] (V3) ..  controls (-0.5,3.5) ..  (Y);
\draw[->] (V2) ..  controls (-1.5,5.6) ..  (Y);
\end{tikzpicture}
  \label{visible:sub2}
\end{subfigure}
\caption{Two configurations where the edge $X \rightarrow Y$ is visible.}
\label{fig:visible}
\end{figure}
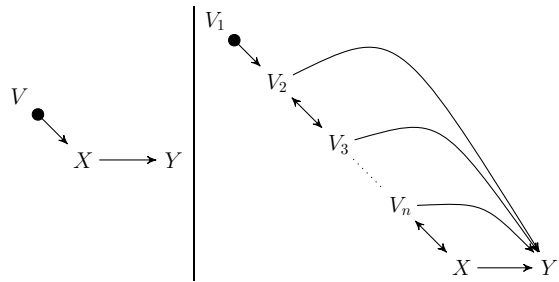

A directed visible edge $X\rightarrow Y$ means that there are no latent confounders between $X$ and $Y$.

\section{MAIN RESULT}\label{sec: main result}

Throughout, let $\g=(\vars,\e)$ represent a $\DAG$, $\CPDAG$, $\MAG$ or $\PAG$, and let $\mathbf{X}$, $\mathbf{Y}$ and $\mathbf{Z}$ be pairwise disjoint subsets of $\vars$, with $\mathbf{X}\neq \emptyset$ and $\mathbf{Y}\neq \emptyset$. Here $\mathbf{X}$ represents the intervention variables and $\mathbf{Y}$ represents the set of response variables, i.e., we are interested in the causal effect of $\mathbf{X}$ on $\mathbf{Y}$. 

We define sound and complete graphical conditions for adjustment sets relative to ($\mathbf{X,Y}$) in $\g$. Thus, if a set $\mathbf{Z}$ satisfies our conditions  relative to ($\mathbf{X},\mathbf{Y}$) in $\g$ (Definition \ref{defacgen}), then it is a valid adjustment set for calculating the causal effect of $\mathbf{X}$ on $\mathbf{Y}$ (Definition \ref{defadjustment}), and every existing valid adjustment set satisfies our conditions (see Theorem \ref{thmgenac}).
First, we define what we mean by an adjustment set.

\begin{definition}{(\textbf{Adjustment set}; \citealp{maathuis2013generalized})}
   Let $\g$ represent a $\DAG$, $\CPDAG$, $\MAG$ or $\PAG$. Then $\mathbf{Z}$ is an adjustment set relative to ($\mathbf{X,Y}$) in $\g$  if for any density $f$ consistent with $\g$ we have
   \resizebox{19.5pc}{!}{
   $
   f(\mathbf{y}|do(\mathbf{x}))=
   \begin{cases}
   f(\mathbf{y}|\mathbf{x}) & \text{if }\mathbf{Z} = \emptyset,\\
   \int_{\mathbf{Z}}f(\mathbf{y}|\mathbf{x,z})f(\mathbf{z})d\mathbf{z} = E_{\mathbf{Z}}\{ f(\mathbf{y}|\mathbf{z,x})\} & \text{otherwise.}
   \end{cases}
   $
   }
   If $\mathbf{X} = \{X\}$ and $\mathbf{Y} = \{Y\}$, we call $\mathbf{Z}$ an adjustment set relative to $(X,Y)$ in the given graph.
   \label{defadjustment}
\end{definition}

To define our generalized adjustment criterion, we introduce the concept of \textit{amenability}:
\begin{definition}{(\textbf{Amenability for $DAG$s, $CPDAG$s, $MAG$s and $PAG$s})}
  A $\DAG$, $\CPDAG$, $\MAG$ or $\PAG$ $\g$ is said to be adjustment amenable, relative to $(\mathbf{X},\mathbf{Y})$ if every possibly directed proper path from $\mathbf{X}$ to $\mathbf{Y}$ in $\g$ starts with a visible edge out of $\mathbf{X}$.
  \label{defadjamengen}
\end{definition}

For conciseness, we will also write ``amenable'' instead of ``adjustment amenable''. The intuition behind the concept of amenability is the following. In $\MAG$s and $\PAG$s, directed edges $X\to Y$ can represent causal effects, but also mixtures of causal effects and latent confounding; in $\CPDAG$s and $\PAG$s, there are edges with unknown direction. This complicates adjustment because paths containing such edges can correspond to causal paths in some represented $\DAG$s and to non-causal paths in others.
For instance, when the graph $X \to Y$ is interpreted as a $\DAG$, the empty set is a valid adjustment set with respect to $(X,Y)$ because there is only one path from $X$ to $Y$, which is causal. When the same graph is however interpreted as a $\MAG$, it can still represent the $\DAG$ $X \to Y$, but also for example the $\DAG$ $X\to Y$ with a non-causal path $X \leftarrow L \to Y$ where $L$ is latent. A similar problem arises in the $\CPDAG$ \begin{tikzpicture}[baseline=-3.2pt,inner sep=0pt,outer sep=2pt]
\node (x) at (0,0) {$X$};
\node (y) at (1,0) {$Y$};
\draw [o-o] (x) -- (y);
\end{tikzpicture}.

We will show that for a graph $\g$ that is not amenable relative to $(\mathbf{X,Y})$, there is no adjustment set relative to $(\mathbf{X,Y})$ in the sense of Definition \ref{defadjustment} (see Lemma \ref{lemmaadjamengen}). Note that every $\DAG$ is amenable, since all edges in a $\DAG$ are visible and directed. For $\MAG$s, our notion of amenability reduces to the one defined by \citet{vanconstructing}.

We now introduce our Generalized Adjustment Criterion (GAC) for $\DAG$s, $CPDAG$s, $\MAG$s and $\PAG$s.

\begin{definition}{(\textbf{Generalized Adjustment Criterion (GAC)})}
   Let $\g$ represent a $\DAG$, $\CPDAG$, $\MAG$ or $\PAG$. Then $\mathbf{Z}$ satisfies the generalized adjustment criterion
      relative to $(\mathbf{X,Y})$ in $\g$ if the following three conditions hold:
   \begin{itemize}
   \item[(0)] $\g$ is adjustment amenable relative to ($\mathbf{X,Y}$), and
   \item[(1)] no element in $\mathbf{Z}$ is a possible descendant in $\g$ of any $W \in \vars\setminus\mathbf{X}$ which lies on a proper possibly causal path from $\mathbf{X}$ to $\mathbf{Y}$, and
   \item[(2)] all proper definite status non-causal paths in $\g$ from $\mathbf{X}$ to $\mathbf{Y}$ are blocked by $\mathbf{Z}$.
   \end{itemize}
   \label{defacgen}
\end{definition}

Note that condition (0) does not depend on $\mathbf{Z}$. In other words, if condition (0) is violated, then there is no set $\mathbf{Z'} \subseteq \vars \setminus (\mathbf{X} \cup\mathbf{Y})$ that satisfies the generalized adjustment criterion relative to $(\mathbf{X,Y})$ in $\g$.

Condition (1) defines a set of nodes that cannot be used in an adjustment set. Denoting this set of forbidden nodes by
\begin{align}\label{eq: def forbidden nodes}
   \fb{\g} = \{ & W' \in \vars: W' \in \PossDe(W,\g) \,\text{for some} \notag \\
   &  W \notin \mathbf{X} \,\text{which lies on a proper possibly} \notag \\
   &  \text{causal path from} \, \mathbf{X} \,\text{to}\, \mathbf{Y}\},
\end{align}
condition (1) can be stated as: $\mathbf{Z} \cap \fb{\g} = \emptyset$.
We will sometimes use this notation in examples and proofs.

We now give the main theorem of this paper.

\begin{theorem}
   Let $\g$ represent a $\DAG$, $\CPDAG$, $\MAG$ or $\PAG$. Then $\mathbf{Z}$ is an adjustment set relative to $(\mathbf{X},\mathbf{Y})$ in $\g$ (Definition \ref{defadjustment}) if and only if $\mathbf{Z}$ satisfies the generalized adjustment criterion  relative to $(\mathbf{X},\mathbf{Y})$ in $\g$ (Definition \ref{defacgen}).
   \label{thmgenac}
\end{theorem}

\section{EXAMPLES}\label{sec: examples}

We now provide some examples that illustrate how the generalized adjustment criterion can be applied.

\begin{example}
   We first return to the example of the Introduction. Consider the $\CPDAG$ $\g[C]$ in Figure~\ref{fig:cpdagexample}a. Note that $\g[C]$ is amenable relative to ($X,Y$) and that $\f{\g[C]} = \{ Y \}$. Hence, any node other than $X$ and $Y$ can be used in an adjustment set. Note that every definite status non-causal path $p$ from $X$ to $Y$ has one of the following paths as a subsequence: $p_1 = \langle X, Z,Y\rangle$ and $p_2 = \langle X, A,B, Y\rangle$, and nodes on $p$ that are not on $p_1$ or $p_2$ are non-colliders on $p$. Hence, if we block $p_1$ and $p_2$, then we block all definite status non-causal paths from $X$ to $Y$. This implies that any superset of $\{Z,A\}$ and $\{Z,B\}$ is an adjustment set relative to $(X,Y)$ in $\g[C]$, and all adjustment sets are given by: $\{Z,A\}$, $\{Z,B\},\{Z,A,I\},\{Z,B,I\},\{Z,A,B\}$ and $\{Z,A,B,I\}$.
  \label{ex1}
\end{example}

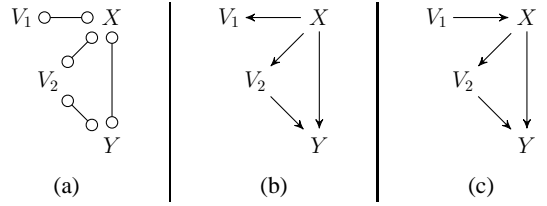
\begin{figure}[!tbp]
   \centering
   \begin{subfigure}{.16\textwidth}
     \centering
     \begin{tikzpicture}[>=stealth',shorten >=1pt,auto,node distance=2cm,main node/.style={minimum size=0.6cm,font=\sffamily\Large\bfseries},scale=0.6,transform shape]
   \node[main node]         (X)                        {$X$};
   \node[main node]         (V1) [left of= X]  		{$V_{1}$};
   \node[main node]         (V2) [below left of = X] 	{$V_{2}$};
   \node[main node]       	 (Y)  [below right of= V2] 	{$Y$};
   \draw[o-o] (V2) edge    (X);
   \draw[o-o] (V1) edge    (X);
   \draw[o-o] (X) edge    (Y);
   \draw[o-o] (V2) edge    (Y);
   \end{tikzpicture}
     \caption{}
     \label{visible:adj1}
   \end{subfigure}%
   \unskip
   \vrule
   \begin{subfigure}{.16\textwidth}
     \centering
     \begin{tikzpicture}[>=stealth',shorten >=1pt,auto,node distance=2cm,main node/.style={minimum size=0.6cm,font=\sffamily\Large\bfseries},scale=.6,transform shape]
   \node[main node]         (X)                        {$X$};
   \node[main node]         (V1) [left of= X]  		{$V_{1}$};
   \node[main node]         (V2) [below left of = X] 	{$V_{2}$};
   \node[main node]       	 (Y)  [below right of= V2] 	{$Y$};
   \draw[<-] (V2) edge    (X);
   \draw[<-] (V1) edge    (X);
   \draw[->] (X) edge    (Y);
   \draw[->] (V2) edge    (Y);
   \end{tikzpicture}
     \caption{}
     \label{visible:adj3}
   \end{subfigure}%
   \unskip
   \vrule
   \begin{subfigure}{.16\textwidth}
     \centering
     \begin{tikzpicture}[>=stealth',shorten >=1pt,auto,node distance=2cm,main node/.style={minimum size=0.6cm,font=\sffamily\Large\bfseries},scale=.6,transform shape]
   \node[main node]         (X)                        {$X$};
   \node[main node]         (V1) [left of= X]  		{$V_{1}$};
   \node[main node]         (V2) [below left of = X] 	{$V_{2}$};
   \node[main node]       	 (Y)  [below right of= V2] 	{$Y$};
   \draw[<-] (V2) edge    (X);
   \draw[->] (V1) edge    (X);
   \draw[->] (X) edge    (Y);
   \draw[->] (V2) edge    (Y);
   \end{tikzpicture}
     \caption{}
     \label{visible:adj2}
   \end{subfigure}
   \caption{(a) $\PAG$ $\g[P]$, (b) $\MAG$ $\g[M]_1$, (c) $\MAG$ $\g[M]_2$ used in Example \ref{ex2}.}
   \label{figex2}
\end{figure}

\begin{example}
   To illustrate the concept of amenability, consider Figure~\ref{figex2} with a $\PAG$ $\g[P]$ in (a), and two $\MAG$s $\g[M]_1$ and $\g[M]_2$ in $[\g[P]]$ in (b) and (c). Note that $\g[P]$ and $\g[M]_1$ are not amenable relative to $(X,Y)$. For $\g[P]$ this is due to the path $X \circcirc Y$, and for $\g[M]_1$ this is due to the invisible edge $X \to Y$. On the other hand, $\g[M]_2$ is amenable relative to $(X,Y)$, since the edges $X\to Y$ and $X \to V_2$ are visible due to the edge $V_1\to X$ and the fact that $V_1$ is not adjacent to $Y$ or $V_2$.
   Since there are no proper definite status non-causal paths from $X$ to $Y$ in $\g[M]_2$, it follows that the empty set satisfies the generalized adjustment criterion relative to $(X,Y)$ in $\g[M]_2$.
   Finally, note that $\g[M]_1$ could also be interpreted as a $\DAG$. In that case it would be amenable relative to $(X,Y)$. This shows that amenability depends crucially on the interpretation of the graph.
   \label{ex2}
\end{example}

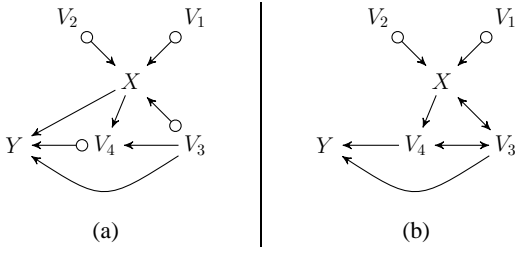
\begin{figure}[!tbp]
   \centering
   \begin{subfigure}{.24\textwidth}
     \centering
     \begin{tikzpicture}[>=stealth',shorten >=1pt,auto,node distance=2cm,main node/.style={minimum size=0.6cm,font=\sffamily\Large\bfseries},scale=0.6,transform shape]
   \node[main node]         (X)                        {$X$};
   \node[main node]         (V2) [above left of= X]  	{$V_{2}$};
   \node[main node]         (V1) [above right of= X]  	{$V_{1}$};
   \node[main node]         (V3) [below right of = X]  			{$V_{3}$};
   \node[main node]         (V4) [left of= V3]  		{$V_{4}$};
   \node[main node]       	 (Y)  [left of= V4]         	{$Y$};
   \draw[o->] (V1) edge    (X);
   \draw[o->] (V2) edge    (X);
   \draw[o->] (V3) edge    (X);
   \draw[->] (V3) edge    (V4);
   \draw[o->] (V4) edge    (Y);
   \draw[->] (X) edge   (Y);
   \draw[->] (X) edge   (V4);
   \draw[->] (V3) ..  controls (-0.6,-2.7) ..   (Y);
   \end{tikzpicture}
     \caption{}
     \label{amen1}
   \end{subfigure}%
   \unskip
   \vrule
   \begin{subfigure}{.24\textwidth}
     \centering
     \begin{tikzpicture}[>=stealth',shorten >=1pt,auto,node distance=2cm,main node/.style={minimum size=0.6cm,font=\sffamily\Large\bfseries},scale=.6,transform shape]
   \node[main node]         (X)                        {$X$};
   \node[main node]         (V2) [above left of= X]  	{$V_{2}$};
   \node[main node]         (V1) [above right of= X]  	{$V_{1}$};
   \node[main node]         (V3) [below right of = X]  {$V_{3}$};
   \node[main node]         (V4) [left of = V3]  		{$V_{4}$};
   \node[main node]       	 (Y)  [left of= V4]         	{$Y$};
   \draw[o->] (V1) edge    (X);
   \draw[o->] (V2) edge    (X);
   \draw[<->] (V3) edge    (X);
   \draw[<->] (V3) edge    (V4);
   \draw[->] (V4) edge    (Y);
   \draw[->] (X) edge   (V4);
   \draw[->] (V3) ..  controls (-0.6,-2.7) ..   (Y);
   \end{tikzpicture}
     \caption{}
     \label{amen2}
   \end{subfigure}
   \caption{(a) $\PAG$ $\g[P]_1$, (b) $\PAG$ $\g[P]_2$ used in Example \ref{ex3}.}
   \label{figex3}
\end{figure}

\begin{example}
  Let $\g[P]_1$ and $\g[P]_2$ be the $\PAG$s in Figure~\ref{figex3}(a) and Figure~\ref{figex3}(b), respectively.
  Both $\PAG$s are amenable relative to $(X,Y)$. We will show that there is an adjustment set relative to $(X,Y)$ in $\g[P]_1$ but not in $\g[P]_2$. This illustrates that amenability is not a sufficient
  criterion for the existence of an adjustment set.

  We first consider $\g[P]_1$. Note that  $\f{\g[P]_1} = \{ V_4, Y \}$ is the set of nodes that cannot be used for adjustment.  There are two proper definite status non-causal paths from $X$ to $Y$ in $\g[P]_1$: $X \arrowcirc V_3 \rightarrow Y$ and $X \rightarrow V_4 \leftarrow V_3 \rightarrow Y$. These are blocked by any set containing $V_3$. Hence, all sets satisfying the GAC relative to ($X,Y$) in $\g[P]_1$ are: $\{ V_3 \}$, $\{V_1, V_3 \}, \{V_2, V_3 \}$ and $\{ V_1, V_2, V_3 \}$.

  We now consider $\g[P]_2$. Note that $\f{\g[P]_2} = \f{\g[P]_1} = \{ V_4, Y \}$. There are three proper definite status non-causal paths from $X$ to $Y$ in $\g[P]_2$:  $p_1 = X \leftrightarrow V_3 \rightarrow Y$, $p_2 = X \leftrightarrow V_3 \leftrightarrow V_4 \to Y$ and  $p_3 = X \rightarrow V_4 \leftrightarrow V_3 \rightarrow Y$. To block $p_1$, we must also use $V_3$. This implies that we must use $V_4$ to block $p_2$. But $V_4 \in \f{\g[P]_2}$. Hence, there is no set $\mathbf{Z}$ that satisfies the GAC relative to ($X,Y$) in $\g[P]_2$.
  \label{ex3}
\end{example}

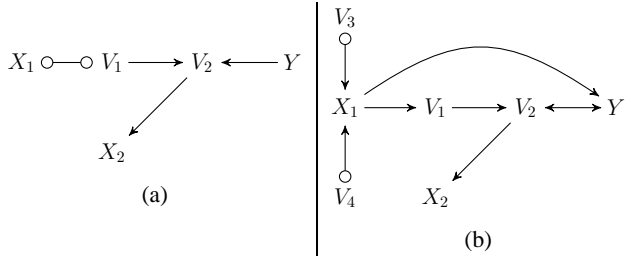
\begin{figure}[!tbp]
   \centering
   \begin{subfigure}{.25\textwidth}
     \centering
     \begin{tikzpicture}[>=stealth',shorten >=1pt,auto,node distance=2cm,main node/.style={minimum size=0.5cm,font=\sffamily\Large\bfseries},scale=0.6,transform shape]
   \node[main node]         (X1)                       {$X_1$};
   \node[main node]         (V1) [right of= X1]  		{$V_{1}$};
   \node[main node]         (X2) [below of= V1]  		{$X_{2}$};
   \node[main node]         (V2) [right of = V1]  		{$V_2$};
   \node[main node]       	 (Y)  [right of= V2]        {$Y$};
   \draw[o-o] (X1) edge    (V1);
   \draw[->] (V1) edge    (V2);
   \draw[->] (Y) edge    (V2);
   \draw[->] (V2) edge    (X2);
   \end{tikzpicture}
     \caption{}
     \label{amenab1}
   \end{subfigure}%
   \unskip
   \vrule
   \begin{subfigure}{.25\textwidth}
     \centering
     \begin{tikzpicture}[>=stealth',shorten >=1pt,auto,node distance=2cm,main node/.style={minimum size=0.5cm,font=\sffamily\Large\bfseries},scale=.6,transform shape]
   \node[main node]         (X1)                       {$X_1$};
   \node[main node]         (V3) [above of= X1]   {$V_3$};
   \node[main node]         (V4) [below of= X1]   {$V_4$};
   \node[main node]         (V1) [right of= X1]  		{$V_{1}$};
   \node[main node]         (V2) [right of= V1]  		{$V_{2}$};
   \node[main node]         (X2) [below of = V1]  		{$X_2$};
   \node[main node]       	 (Y)  [right of= V2]       	{$Y$};
   \draw[->] (X1) edge    (V1);
   \draw[o->] (V3) edge    (X1);
   \draw[o->] (V4) edge    (X1);
   \draw[->] (V1) edge    (V2);
   \draw[<->] (Y) edge    (V2);
   \draw[->] (V2) edge    (X2);
   \draw[->] (X1) ..  controls (2.5,1.7) and (3.5,1.7) ..   (Y);
   \end{tikzpicture}
     \caption{}
     \label{amenab2}
   \end{subfigure}
   \caption{(a) $\CPDAG$ $\g[C]$, (b) $\PAG$ $\g[P]$ used in Example \ref{ex4}.}
   \label{figex4}
\end{figure}

\begin{example}
 Let $\mathbf{X} = \{X_1,X_2\}$ and $\mathbf{Y} = \{Y\}$ and consider the $\CPDAG$ $\g[C]$ and the $\PAG$ $\g[P]$ in Figures \ref{figex4}(a) and \ref{figex4}(b). We will show that for both graphs there is no set that satisfies the generalized back-door criterion of  \cite{maathuis2013generalized} relative to $(\mathbf{X},\mathbf{Y})$, but there are sets that satisfy the generalized adjustment criterion relative to $(\mathbf{X},\mathbf{Y})$ in these graphs.
  
   Recall that a set $\mathbf{Z}$ satisfies the generalized back-door criterion relative to $(\mathbf{X},\mathbf{Y})$ and a $\CPDAG$ ($\PAG$) $\g$ if $\mathbf{Z}$ contains no possible descendants of $\mathbf{X}$ in $\g$ and if for every $X \in \mathbf{X}$ the set $\mathbf{Z}\cup \mathbf{X}\setminus\{X\}$ blocks every definite status path from $X$ to every $Y \in \mathbf{Y}$ in $\g$ that does not start with a visible edge out of $X$.

   We first consider the $\CPDAG$ $\g[C]$. 
   To block the path $X_2 \leftarrow V_2 \leftarrow Y$, we must use node $V_2$, but $V_2 \in \PossDe(X_1, \g[C])$. Hence, no set $\mathbf{Z}$ can satisfy the generalized back-door criterion relative to $(\mathbf{X},\mathbf{Y})$ in $\g[C]$. However, $\{ V_1,V_2\}$ satisfies the generalized adjustment criterion relative to $(\mathbf{X},\mathbf{Y})$ in $\g[C]$.

We now consider $\g[P]$. 
To block the path $X_2 \leftarrow V_2 \leftrightarrow Y$, we must use node $V_2$. But, $V_2 \in \De(X_1,\g[P])$ and thus there is no set satisfying the generalized back-door criterion relative to $(\mathbf{X},\mathbf{Y})$ in $\g[P]$. However, sets $\{ V_1,V_2\}$, $\{V_1,V_2,V_3\}$, $\{V_1,V_2,V_4\}$, $\{V_1,V_2,V_3,V_4\}$ all satisfy the generalized adjustment criterion relative to $(\mathbf{X},\mathbf{Y})$ in $\g[P]$. 
   \label{ex4}
\end{example}

\section{PROOF OF THEOREM \ref{thmgenac}}\label{sec: proof}

For $\DAG$s and $\MAG$s, our generalized adjustment criterion reduces to the following adjustment criterion:
\begin{definition}{(\textbf{Adjustment Criterion (AC)})}
   Let $\g = (\vars,\e)$ represent a $\DAG$ or $\MAG$ . Then $\mathbf{Z}$ satisfies the adjustment criterion
    relative to $(\mathbf{X},\mathbf{Y})$ in $\g$ if the following three conditions hold:
   \begin{itemize}
   \item[(0*)] $\g$ is adjustment amenable with respect to $(\mathbf{X},\mathbf{Y})$, and
   \item[(a)] no element in $\mathbf{Z}$ is a descendant in $\g$ of any $W \in \vars \setminus \mathbf{X}$ which lies on a proper causal path from $\mathbf{X}$ to $\mathbf{Y}$, and
   \item[(b)] all proper non-causal paths in $\g$ from $\mathbf{X}$ to $\mathbf{Y}$ are blocked by $\mathbf{Z}$.
   \end{itemize}
   \label{defac}
\end{definition}
This adjustment criterion is a slightly reformulated but equivalent version of the adjustment criterion of \citet{shpitser2012validity} for $\DAG$s and of \citet{vanconstructing} for $\MAG$s, with amenability directly included in the criterion. This adjustment criterion was shown to be sound and complete for $\DAG$s \citep{shpitser2012validity,shpitser2012avalidity} and $\MAG$s \citep{vanconstructing}. We therefore only need to prove Theorem \ref{thmgenac} for $\CPDAG$s and $\PAG$s.

To this end, we need three main lemmas, given below.
Throughout, we let $\mathcal G=(\vars, \e)$ represent a $\CPDAG$ or a $\PAG$, and we let $\mathbf{X}$, $\mathbf{Y}$ and $\mathbf{Z}$ be pairwise disjoint subsets of $\vars$, with $\mathbf{X}\neq \emptyset$ and $\mathbf{Y}\neq \emptyset$. We use GAC and AC to refer to the generalized adjustment criterion (Definition \ref{defacgen}) and adjustment criterion (Definition \ref{defac}), respectively.

Lemma \ref{lemmaadjamengen} is about condition (0) of the GAC:
\begin{lemma}
    If a $\CPDAG$ ($\PAG$) $\g$ satisfies condition (0) of the GAC relative to $(\mathbf{X},\mathbf{Y})$, then every $\DAG$ ($\MAG$) in $[\g]$ satisfies condition (0*) of the AC relative to $(\mathbf{X},\mathbf{Y})$. On the other hand, if $\g$ violates condition (0) of the GAC relative to $(\mathbf{X},\mathbf{Y})$, then there exists no set $\mathbf{Z}' \subseteq \vars \setminus (\mathbf{X}\cup \mathbf{Y})$ that is an adjustment set relative to $(\mathbf{X},\mathbf{Y})$ in $\g$ (see Definition \ref{defadjustment}).
   \label{lemmaadjamengen}
\end{lemma}

Next, we assume that $\g$ satisfies condition (0) of the GAC relative to $(\mathbf{X},\mathbf{Y})$.  Under this assumption, we show that $\mathbf{Z}$ satisfies conditions (1) and (2) of the GAC relative to $(\mathbf{X,Y})$ in $\g$ if and only if $\mathbf{Z}$ satisfies conditions (a) and (b) of the AC relative to $(\mathbf{X,Y})$ in every $\DAG$ ($\MAG$) in $[\g]$. This is shown in two separate lemmas:

\begin{lemma}
   Let condition (0) of the GAC be satisfied relative to $(\mathbf{X},\mathbf{Y})$ in a $\CPDAG$ ($\PAG$) $\g$. Then the following two statements are equivalent:
   \begin{itemize}
       \item $\mathbf{Z}$ satisfies condition (1) of the GAC relative to $(\mathbf{X},\mathbf{Y})$ in $\g$.
       \item $\mathbf{Z}$ satisfies condition (a) of the AC relative to $(\mathbf{X},\mathbf{Y})$ in every $\DAG$ ($\MAG$) in $[\g]$.
   \end{itemize}
   \label{lemmaprovingeqofacconditiona}
\end{lemma}

\begin{lemma}
   Let condition (0) of the GAC be satisfied relative to $(\mathbf{X},\mathbf{Y})$ in a $\CPDAG$ ($\PAG$) $\g$, and let $\mathbf{Z}$ satisfy condition (1) of the GAC relative to  $(\mathbf{X},\mathbf{Y})$ in $\g$. Then the following two statements are equivalent:
   \begin{itemize}
      \item $\mathbf{Z}$ satisfies condition (2) of the GAC relative to $(\mathbf{X},\mathbf{Y})$ in $\g$.
      \item $\mathbf{Z}$ satisfies condition (b) of the AC relative to $(\mathbf{X},\mathbf{Y})$ in every $\DAG$ ($\MAG$) in $[\g]$.
   \end{itemize}
   \label{lemmaequivalenceofcondb}
\end{lemma}

The proofs of Lemmas \ref{lemmaadjamengen}, \ref{lemmaprovingeqofacconditiona} and \ref{lemmaequivalenceofcondb} are discussed in Sections \ref{sec: adjustment amenable proof}, \ref{sec: eq of condition a} and \ref{sec: eq of condition b}, respectively. Some proofs require additional lemmas that can be found in the supplement. The proof of Lemma \ref{lemmaequivalenceofcondb} is the most technical, and builds on the work of \citet{zhang2006causal}.

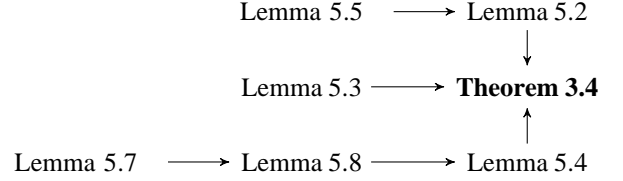
\begin{figure}[!tbp]
    \centering
      \begin{tikzpicture}[>=stealth',shorten >=1pt,node distance=3cm, main node/.style={minimum size=0.4cm}]
    \node[main node]         (T34) {\textbf{Theorem \ref{thmgenac}}};
    \node[main node,yshift=1cm] (L52) at (T34) {Lemma \ref{lemmaadjamengen}};
    \node[main node]         (L31)  [left of= L52,text width=2.2cm,align=center]          {Lemma \ref{lemmadirectedpath}};
    \node[main node]         (L53)  [left of= T34]          {Lemma \ref{lemmaprovingeqofacconditiona}};
    \node[main node,yshift=-1cm] (L54) at (T34) {Lemma \ref{lemmaequivalenceofcondb}};
    \node[main node]            (L56)  [left of= L54]   {Lemma \ref{lemma22}};
    \node[main node]            (L57)  [left of= L56,text width=2.2cm,align=center]   {Lemma \ref{lemma11}};
    \draw[->] (L52) edge    (T34);
    \draw[->] (L53) edge    (T34);
    \draw[->] (L31) edge    (L52);
    \draw[->] (L54) edge    (T34);
    \draw[->] (L56) edge    (L54);
    \draw[->] (L57) edge    (L56);
    \end{tikzpicture}
    \caption{Proof structure of Theorem \ref{thmgenac}.}
    \label{figproof}
\end{figure}

Figure~\ref{figproof} shows how all lemmas fit together to prove Theorem \ref{thmgenac}.

\begin{proofthm}
   First, suppose that the $\CPDAG$ ($\PAG$) $\g$ and the sets $\mathbf{X}$, $\mathbf{Y}$ and $\mathbf{Z}$ satisfy all conditions of the GAC. By applying Lemmas \ref{lemmaadjamengen}, \ref{lemmaprovingeqofacconditiona} and \ref{lemmaequivalenceofcondb} in turn, it directly follows that all conditions of the AC are satisfied by $\mathbf{X}$, $\mathbf{Y}$ and $\mathbf{Z}$ and any $\DAG$ ($\MAG$) in $[\g]$.

   To prove the other direction, suppose that the tuple $\g$, $\mathbf{X}$, $\mathbf{Y}$, $\mathbf{Z}$ does not satisfy all conditions of the GAC. First, suppose that $\g$ violates condition (0) relative to $(\mathbf{X},\mathbf{Y})$. Then by Lemma \ref{lemmaadjamengen}, there is no adjustment set relative to $(\mathbf{X},\mathbf{Y})$ in $\g$, and hence $\mathbf{Z}$ is certainly not an adjustment set.

   Otherwise, $\mathbf{Z}$ must violate condition (1) or (2) of the GAC relative to $(\mathbf{X},\mathbf{Y})$. By applying Lemmas \ref{lemmaprovingeqofacconditiona} and \ref{lemmaequivalenceofcondb} in turn, this implies that there is a $\DAG$ $\g[D]$ ($\MAG$ $\g[M]$) in $[\g]$ such that $\mathbf{Z}$ violates conditions (a) or (b) of the AC relative to $(\mathbf{X},\mathbf{Y})$ in $\g[D]$ ($\g[M]$).
   Since the AC is sound and complete for $\DAG$s and $\MAG$s, this implies that $\mathbf{Z}$ is not an adjustment set relative to $(\mathbf{X},\mathbf{Y})$ in $\g[D]$ ($\g[M]$), so that $\mathbf{Z}$ is certainly not an adjustment set relative to $(\mathbf{X},\mathbf{Y})$ in $\g$.

\end{proofthm}

\subsection{PROOF OF LEMMA \ref{lemmaadjamengen}}\label{sec: adjustment amenable proof}

The proof of Lemma 5.2 is based on the following lemma:

\begin{lemma}
 Let $X$ and $Y$ be nodes in a $\PAG$ $\g[P]$, such that there is a possibly directed path $\pstar$ from $X$ to $Y$ in $\g[P]$ that does not start with a visible edge out of $X$. Then there is a $\MAG$ $\g[M]$ in $[\g[P]]$ such that the path $p$ in $\g[M]$, consisting of the same sequence of nodes as $\pstar$ in $\g[P]$, contains a subsequence that is a directed path from $X$ to $Y$ starting with an invisible edge in $\g[M]$.
    \label{lemmadirectedpath}
\end{lemma}

The proof of Lemma \ref{lemmadirectedpath} is given in the supplement.

\begin{proofsketchamenab}
   First suppose that $\g$ satisfies condition (0) of the GAC relative to $(\mathbf{X},\mathbf{Y})$, meaning that every proper possibly directed path from $\mathbf{X}$ to $\mathbf{Y}$ in $\g$ starts with a visible edge out of $\mathbf{X}$. Any visible edge in $\g$ is visible in all $\DAG$s ($\MAG$s) in $[\g]$, and any proper directed path in a $\DAG$ ($\MAG$) in $[\g]$ corresponds to a proper possibly directed path in $\g$. Hence, any proper directed path from $\mathbf{X}$ to $\mathbf{Y}$ in any $\DAG$ ($\MAG$) in $[\g]$ starts with a visible edge out of $\mathbf{X}$. This shows that all $\DAG$s ($\MAG$s) in $[\g]$ satisfy condition (0*) of the AC relative to $(\mathbf{X},\mathbf{Y})$.

   Next, suppose that $\g$ violates condition (0) of the GAC relative to $(\mathbf{X},\mathbf{Y})$. We will show that this implies that there is no set $\mathbf{Z}' \subseteq  \vars \setminus (\mathbf{X} \cup \mathbf{Y})$ that is an adjustment set relative to $(\mathbf{X},\mathbf{Y})$ in $\g$.
   We give separate proofs for $\CPDAG$s and $\PAG$s.

   Thus, let $\g$ represent a $\CPDAG$ and suppose that there is a proper possibly directed path $p$ from a node $X \in \mathbf{X}$ to a node $Y \in \mathbf{Y}$ that starts with a non-directed edge ($\circcirc$).

   Let $p' = \langle X, V_1,\dots, Y \rangle$ (where $V_1 = Y$ is allowed) be a shortest subsequence of $p$ such that $p'$ is also a proper possibly directed path from $X$ to $Y$ starting with a non-directed edge in $\g$. We first show that $p'$ is a definite status path, by contradiction. Thus, suppose that $p'$ is not a definite status path. Then the length of $p'$ is at least 2, and we write $p' = \langle X,V_1,\dots,V_k=Y\rangle$ for $k\ge 2$. Since the subpath $p'(V_1,Y)$ is a definite status path (otherwise we can choose a shorter path), this means that $V_1$ is not of a definite status on $p'$. This implies the existence of an edge between $X$ and $V_2$. This edge must be of the form $X \to V_2$, since $X \circcirc V_2$ implies that we can choose a shorter path, and $X \leftarrow V_2$ together with $X \circcirc V_1$ implies $V_1 \leftarrow V_2$  by Lemma 1 from \citet{meek1995causal} (see Section \ref{S-sec:1} of the supplement), so that $p'$ is not possibly directed from $X$ to $Y$. But the edge $X \to V_2$ implies that $V_1 \to V_2$, since otherwise Lemma 1 from \citet{meek1995causal} implies $X \to V_1$. But then $V_1$ is a definite non-collider on $p'$, which contradicts that $V_1$ is not of definite status.

   Hence, $p'$ is a proper possibly directed definite status path from $X$ to $Y$. By Lemma 7.6 from \citet{maathuis2013generalized} (see Section \ref{S-sec:1} of the supplement), there is a $\DAG$ $\g[D]_1$ in $[\g]$ such that there are no additional arrowheads into $X$, as well as a $\DAG$ $\g[D]_2$ in $[\g]$ such that there are no additional arrowheads into $V_1$. This means that the paths corresponding to $p'$ are oriented as $p'_1 = X \to V_1 \to \dots \to Y$ and $p'_2 = X \leftarrow V_1 \to \dots \to Y$ in $\g[D]_1$ and $\g[D]_2$.
   An adjustment set relative to $(\mathbf{X},\mathbf{Y})$ in $\g[D]_2$ must block the non-causal path $p'_{2}$, by using at least one of the non-endpoints nodes on this path. But all these nodes are in $\fb{\g[D]_1}$ (see \eqref{eq: def forbidden nodes}). Hence, there is no set $\mathbf{Z}' \subseteq \vars \setminus (\mathbf{X} \cup \mathbf{Y})$ that satisfies the AC relative to $(\mathbf{X},\mathbf{Y})$ in $\g[D]_1$ and $\g[D]_2$ simultaneously. Since the AC is sound and complete for $\DAG$s, this implies that there is no $\mathbf{Z}' \subseteq \vars \setminus (\mathbf{X} \cup \mathbf{Y})$ that is an adjustment set relative to $(\mathbf{X},\mathbf{Y})$ in $\g$.

   Finally, let $\g$ represent a $\PAG$ and suppose that there is a proper possibly directed path $p$ from some $X\ \in \mathbf{X}$ to some $Y\in \mathbf{Y}$ that does not start with a visible edge out of $X$ in $\g$.

  By Lemma \ref{lemmadirectedpath}, there is a subsequence $p'$ of $p$ such that there is a $\MAG$ $\g[M]$ in $[\g]$ where the corresponding path is directed from $X$ to $Y$ and starts with an invisible edge.
   Then $\g[M]$ is not amenable relative to $(\mathbf{X},\mathbf{Y})$. By Lemma 5.7 from \citet{vanconstructing} (see Section \ref{S-sec:1} of the supplement) this means that there is no set $\mathbf{Z}' \subseteq \vars \setminus (\mathbf{X} \cup \mathbf{Y})$ that is an adjustment set relative to $(\mathbf{X},\mathbf{Y})$ in $\g[M]$. Hence, there is no set $\mathbf{Z}' \subseteq \vars \setminus (\mathbf{X} \cup \mathbf{Y})$ that is an adjustment set relative to $(\mathbf{X},\mathbf{Y})$ in $\g$.
\end{proofsketchamenab}

\subsection{PROOF OF LEMMA \ref{lemmaprovingeqofacconditiona}}\label{sec: eq of condition a}

\begin{prooflemmaeqa}
   First, suppose that $\mathbf{Z}$ satisfies condition (1) of the GAC relative to $(\mathbf{X},\mathbf{Y})$ in $\g$. Then $\mathbf{Z} \cap \fb{\g} = \emptyset$.
   Since $\fb{\g[D]} \subseteq \fb{\g}$ ($\fb{\g[M]} \subseteq \fb{\g}$) for any $\DAG$ $\g[D]$ ($\MAG$ $\g[M]$) in $[\g]$, it follows directly that $\mathbf{Z}$ satisfies condition (a) of the AC relative to $(\mathbf{X},\mathbf{Y})$ in all $\DAG$s ($\MAG$s) in $[\g]$.

   To prove the other direction, suppose that $\g$ satisfies condition (0) of the GAC relative to $(\mathbf{X},\mathbf{Y})$, but that $\mathbf{Z}$ does not satisfy condition (1) of the GAC relative to $(\mathbf{X,Y})$ in $\g$. Then there is a node $V \in \mathbf{Z} \cap \fb{\g}$, i.e., $V\in \mathbf{Z}$ and $V$ is a possible descendant of a node $W$ on a proper possibly directed path from some $X\in \mathbf{X}$ to some $Y\in \mathbf{Y}$ in $\g$. We denote this path by $p = \langle X, V_1,\dots, V_k, Y \rangle$, where $k\ge 1$ and $W\in \{V_1,\dots,V_k\}$. Then the subpaths $q=p(X,W)$ and $r=p(W,Y)$ are also proper possibly directed paths. Moreover, there is a possibly directed path $s$ from $W$ to $V$, where this path is allowed to be of zero length (if $W=V$). We will show that the existence of these paths implies that there is a $\DAG$ $\g[D]$ ($\MAG$ $\g[M]$) in $[\g]$ such that $\mathbf{Z}$ violates condition (a) of the AC relative to $(\mathbf{X},\mathbf{Y})$ in $\g[D]$ ($\g[M]$).

   By Lemma B.1 from \citet{zhang2008completeness} (see Section \ref{S-sec:1} of the supplement), there are subsequences $q'$, $r'$ and $s'$ of $q$, $r$ and $s$ that are unshielded proper possibly directed paths (again $s'$ is allowed to be a path of zero length). Moreover, $q'$ must start with a directed (visible) edge, since otherwise the concatenated path $q' \oplus r'$, which is again a proper possibly directed path from $X$ to $Y$, would violate condition (0) of the GAC.

 Lemma B.1 from \citet{zhang2008completeness} then implies that $q'$ is a directed path from $X$ to $W$ in $\g$. Hence, the path corresponding to $q'$ is a directed path from $X$ to $W$ in any $\DAG$ ($\MAG$) in $[\g]$.

   By Lemma 7.6 from \citet{maathuis2013generalized}, there is at least one $\DAG$ $\g[D]$ ($\MAG$ $\g[M]$) in  $[\g]$ that has no additional arrowheads into $W$. In this graph $\g[D]$ ($\g[M]$), the path corresponding to $r'$ is a directed path from $W$ to $Y$, and the path corresponding to $s'$ is a directed path $W$ to $V$.
   Hence,  $V \in \fb{\g[D]}$  ($V \in \fb{\g[M]}$), so that $\mathbf{Z}$ does not satisfy condition (a) of the AC relative to $(\mathbf{X},\mathbf{Y})$ in $\g[D]$ ($\g[M]$).
\end{prooflemmaeqa}

\subsection{PROOF OF LEMMA \ref{lemmaequivalenceofcondb}} \label{sec: eq of condition b}

We first define a distance between a path and a set in Definition \ref{defdistancefromC}. We then give the proof of Lemma \ref{lemmaequivalenceofcondb}. This proof relies on Lemma \ref{lemma11} and Lemma \ref{lemma22} which are given later in this section.

\begin{definition}
   (\textbf{Distance-from-$\mathbf{Z}$}; \citealp{zhang2006causal}) Given a path $p$ from $\mathbf{X}$ to $\mathbf{Y}$ that is m-connecting given $\mathbf{Z}$ in a $\DAG$ or $\MAG$, for every collider $Q$ on $p$, there is a directed path (possibly of zero length) from $Q$ to a member of $\mathbf{Z}$. Define the $\distancefrom{\mathbf{Z}}$ of $Q$ to be the length of a shortest directed path (possibly of length 0) from $Q$ to $\mathbf{Z}$, and define the $\distancefrom{\mathbf{Z}}$ of $p$ to be the sum of the distances from $\mathbf{Z}$ of the colliders on $p$.
   \label{defdistancefromC}
\end{definition}

\begin{prooflemmaeqb}
   Let $\g$ represent an amenable $\CPDAG$ ($\PAG$) that satisfies condition (0) of the GAC relative to $(\mathbf{X},\mathbf{Y})$, and let $\mathbf{Z}$ satisfy condition (1) of the GAC relative to $(\mathbf{X},\mathbf{Y})$ in $\g$.

   We first prove that if $\mathbf{Z}$ does not satisfy condition (2) of the GAC relative to $(\mathbf{X},\mathbf{Y})$ in $\g$, then $\mathbf{Z}$ does not satisfy condition (b) of the AC relative to $(\mathbf{X},\mathbf{Y})$ in any $\DAG$ ($\MAG$) in $[\g]$.
   Thus, assume that there is a proper definite status non-causal path $p$ from $X \in \mathbf{X}$ to $Y \in \mathbf{Y}$ that is m-connecting given $\mathbf{Z}$ in $\g$. Consider any $\DAG$ $\g[D]$ ($\MAG$ $\g[M]$) in $[\g]$. Then the path corresponding to $p$ in $\g[D]$ ($\g[M]$) is a proper non-causal m-connecting path from $\mathbf{X}$ to $\mathbf{Y}$ given $\mathbf{Z}$. Hence,
   $\mathbf{Z}$ violates condition (b) of the AC relative to $(\mathbf{X},\mathbf{Y})$ and $\g[D]$ ($\g[M]$).

   Next, we prove that if $\mathbf{Z}$ violates condition (b) of the AC relative to $(\mathbf{X},\mathbf{Y})$ in some $\DAG$ ($\MAG$) in $[\g]$, then $\mathbf{Z}$ violates condition (2) of the GAC relative to $(\mathbf{X},\mathbf{Y})$ in $\g$.
   Thus, assume that there is a $\DAG$ $\g[D]$ ($\MAG$ $\g[M]$) in $[\g]$ such that there is a proper non-causal m-connecting path from $\mathbf{X}$ to $\mathbf{Y}$ in $\g[D]$ ($\g[M]$) given $\mathbf{Z}$.
   We choose a shortest such path $p$, such that no equally short proper non-causal m-connecting path has a shorter $\distancefrom{\mathbf{Z}}$ than $p$. By Lemma \ref{lemma22} below, the corresponding path $\pstar$ in $\g$ is an m-connecting proper definite status non-causal path from $\mathbf{X}$ to $\mathbf{Y}$ given $\mathbf{Z}$. Hence $\mathbf{Z}$ violates condition (b) of the GAC relative to $(\mathbf{X},\mathbf{Y})$ in $\g$.
\end{prooflemmaeqb}

\begin{lemma}
   Let $\g[M]$ represent a $\MAG$ ($\DAG$) and let $\g[P]$ be the $\PAG$ ($\CPDAG$) of $\g[M]$. Let $\g[P]$ satisfy condition (0) of the GAC relative to $(\mathbf{X},\mathbf{Y})$,  and let $\mathbf{Z}$ satisfy condition (1) of the GAC relative to $(\mathbf{X},\mathbf{Y})$ in $\g[P]$. Let $p$ be a shortest proper non-causal path from $\mathbf{X}$ to $\mathbf{Y}$ that is m-connecting given $\mathbf{Z}$ in $\g[M]$ and let $\pstar$ denote the corresponding path constituted by the same sequence of variables in $\g[P]$. Then $\pstar$ is a proper definite status non-causal path in $\g[P]$.
   \label{lemma11}
\end{lemma}

 Lemma \ref{lemma11} is related to Lemma 1 from \citet{zhang2006causal}. The proof of Lemma \ref{lemma11} is given in the supplement.

\begin{lemma}
   Let $\g[M]$ represent a $\MAG$ ($\DAG$) and let $\g[P]$ be the $\PAG$ ($\CPDAG$) of $\g[M]$. Let $\g[P]$ satisfy condition (0) of the GAC relative to $(\mathbf{X},\mathbf{Y})$, and let $\mathbf{Z}$ satisfy condition (1) of the GAC relative to $(\mathbf{X},\mathbf{Y})$ in $\g[P]$. Let $p$ be a shortest proper non-causal path from $\mathbf{X}$ to $\mathbf{Y}$ that is m-connecting given $\mathbf{Z}$ in $\g[M]$, such that no equally short such path has a shorter $\distancefrom{\mathbf{Z}}$ than $p$. Let $\pstar$ denote the corresponding path constituted by the same sequence of variables in $\g[P]$. Then $\pstar$ is a proper definite status non-causal path from $\mathbf{X}$ to $\mathbf{Y}$ that is m-connecting given $\mathbf{Z}$ in $\g[P]$. 
   \label{lemma22}
\end{lemma}

Lemma \ref{lemma22} is is related to Lemma 2 from \citet{zhang2006causal}.

  \begin{proofof}[Lemma \ref{lemma22}]
By Lemma \ref{lemma11}, $\pstar$ is a proper definite status non-causal path in $\g[P]$. It is only left to prove that $\pstar$ is m-connecting given $\mathbf{Z}$ in $\g[P]$.
   
Every definite non-collider on $\pstar$ in $\g[P]$ corresponds to a non-collider on $p$ in $\g[M]$, and every collider on $\pstar$ is also a collider on $p$. Since $p$ is m-connecting given $\mathbf{Z}$, no non-collider is in $\mathbf{Z}$ and every collider has a descendant in $\mathbf{Z}$. Let $Q$ be an arbitrary collider (if there is one). Then there is a directed path (possibly of zero length) from $Q$ to a node in $\mathbf{Z}$ in $\g[M]$. Let $d$ be a shortest such path from $Q$ to a node $Z \in \mathbf{Z}$. Let $\pstar[d]$ denote the corresponding path in $\g[P]$, constituted by the same sequence of variables. Then $\pstar[d]$ is an unshielded possibly directed path from $Q$ to $Z$ in $\g[P]$ (Lemma B.1 from \citet{zhang2008completeness}).
   
   It is only left to show that $\pstar[d]$ is a directed path. If $\pstar[d]$ is of zero length, this is trivially true. Otherwise,
  suppose for contradiction that there is a circle mark on $\pstar[d]$. Then $\pstar[d]$ must start with a circle mark at $Q$ 
  (cf.\ Lemma B.2 from \citealp{zhang2008completeness} and Lemma 7.2 from \citealp{maathuis2013generalized}; see Section \ref{S-sec:1} of the supplement).
  
 Let $S$ be the first node on $d$ after $Q$.
 If $S$ is not a node on $p$, then following the proof of Lemma 2 from \citet{zhang2006causal} there is a path  $p' = p(X,W) \oplus W \bulletarrow S \arrowbullet V \oplus p(V,Y)$, where $W$ and $V$ are nodes distinct from $Q$ on $p(X,Q)$ and $p(Q,Y)$ respectively and $p'$ is m-connecting given $\mathbf{Z}$ in $\g[M]$.
 Since $p'$ is non-causal and shorter than $p$, or as long as $p$ but with a shorter $\distancefrom{\mathbf{Z}}$ than $p$, the path $p'$ must be non-proper, i.e. $S \in \mathbf{X}$.  But, in that case the path $\langle S,V \rangle \oplus p(V,Y)$ is a proper non-causal m-connecting path from $\mathbf{X}$ to $\mathbf{Y}$ given $\mathbf{Z}$ that is shorter than $p$ in $\g[M]$. This contradicts our assumption about $p$.

 If $S$ is a node on $p$, then it lies either on $p(X,Q)$ or $p(Q,Y)$. Assume without loss of generality that $S$ is on $p(Q,Y)$. Following the proof of Lemma 2 from \citet{zhang2006causal}, there exists a path $p' =p(X,W) \oplus  W  \bulletarrow S  \oplus p(S,Y)$ in $\g[M]$, where $W$ is a node on $p(X,Q)$ distinct from $Q$ that is m-connecting given $\mathbf{Z}$ in $\g[M]$.
Since $p'$  is proper, and shorter than $p$, or as long as $p$ but with a shorter $\distancefrom{\mathbf{Z}}$ than $p$, the path $p'$ must be causal in $\g[M]$. Let $\pstar[p']$ denote the corresponding path constituted by the same sequence of variables in $\g[P]$.  Then $\pstar[p']$ is a possibly causal path and $Z \in \PossDe(S,\g[P])$, so $Z \in \fb{\g[P]} \cap \mathbf{Z}$. This is in contradiction with our assumption of $\mathbf{Z}$ satisfying condition ($1$) of the GAC relative to ($\mathbf{X,Y}$) in $\g[P]$.

 Thus, the path $\pstar[d]$ is directed and $Q$ is an ancestor of $\mathbf{Z}$ in $\g[P]$. This proves that $\pstar$ is a proper definite status non-causal path from $\mathbf{X}$ to $\mathbf{Y}$ that is m-connecting given $\mathbf{Z}$ in $\g[M]$.
\end{proofof}

\section{DISCUSSION}

We have derived a generalized adjustment criterion that is necessary and sufficient for adjustment in $\DAG$s,
$\MAG$s, $\CPDAG$s and $\PAG$s. Our criterion unifies existing criteria for $\DAG$s and $\MAG$s, and provides a new result
for $\CPDAG$s and $\PAG$s, where only a sufficient criterion existed until now. This is relevant in practice,
in particular in combination with algorithms that can learn $\CPDAG$s or $\PAG$s from observational data.

Our generalized adjustment criterion is stated in terms of paths that need
to be blocked, which is intuitively appealing. A
logical next step for future research would be to transform
our criterion into an algorithmically constructive version that
could be used to efficiently perform tasks like enumeration
of all minimal adjustment sets for a given graph. This has already
been done for $\DAG$s and $\MAG$s by \citet{vanconstructing},
and we strongly suspect that their results can be extended
to $\CPDAG$s and $\PAG$s as well. In a similar spirit, it would be desirable to
have an easily checkable condition to determine if there exists any adjustment
set at all, as done for the generalized back-door criterion for single interventions by \citet{maathuis2013generalized}.
In turn, these results could then be used to characterize distances between graphs, as done by \citet{Peters2014neco}.
Future work might also explore
under which circumstances our restriction
to not allow for latent selection variables might be relaxed,
or whether our criterion could be combined with methods to
recover from selection bias \citep{Bareinboim2014}.

As pointed out in Section~\ref{sec: examples}, our criterion sometimes has to interpret
$\PAG$s or $\MAG$s differently than $\DAG$s or $\CPDAG$s.
This is the case precisely when the first edge on some proper
possibly
causal path in a $\MAG$ or $\PAG$ is not visible.
However, this difference in interpretation
is irrelevant for $\DAG$s or $\CPDAG$s that would be
amenable when viewed as a $\MAG$ or $\PAG$.
For instance, if we are given a $\DAG$ $\g[D]$
that is amenable when interpreted as a $\MAG$ $\g[M]$,
then its
adjustment sets also work for every $\DAG$ that the $\MAG$ $\g[M]$
represents, many of which could contain
latent confounding variables. Reading a $\DAG$ as a $\MAG$
(or a $\CPDAG$ as a $\PAG$) can thus allow computing adjustment
sets that are to some extent invariant to confounding.

We note that an adjustment set relative to $(\mathbf{X},\mathbf{Y})$ in a given graph can only exist if the total causal effect of $\mathbf{X}$ on $\mathbf{Y}$ is identifiable in the graph. If the effect of $\mathbf{X}$ on $\mathbf{Y}$ is  not identifiable, one may be interested in computing all possible total causal effects of $\mathbf{X}$ on $\mathbf{Y}$ for $\DAG$s represented by the given graph. Such an approach is used in the IDA algorithm of \cite{MaathuisKalischBuehlmann09,MaathuisColomboKalischBuehlmann10}, by  considering all $\DAG$s represented by a $\CPDAG$ and applying back-door adjustment to each of these $\DAG$s. Similar ideas could be used for $\MAG$s and $\PAG$s, but listing all relevant $\DAG$s described by a $\MAG$ or $\PAG$ seems rather non-trivial.

There is also an interesting connection between amenability and instrumental variables:
a $\MAG$ or $\PAG$ $\g$ with $\mathbf{X}=\{X\}$
is amenable with respect to $(\mathbf{X},\mathbf{Y})$ whenever it
contains an \emph{instrument} $I$,
i.e.\ there exists a variable that is a parent of $X$ but not a parent
of any child of $X$ (e.g., $I$ in Figure~\ref{fig:cpdagexample}a).
Thus, instruments are useful to find
adjustment sets in nonparametric graphical models
that allow for latent confounding.
This connection is perhaps surprising
given that the notion of instruments originates from causal effect
identifications in linear models \citep{AngristIR96}.

In summary, our generalized adjustment criterion exhaustively
characterizes the options
to identify total causal effects by covariate adjustment
in $\DAG$s, $\MAG$s, $\CPDAG$s, and $\PAG$s.
Our results entail several existing, less general or less powerful ones
\citep{pearl1993bayesian,shpitser2012validity,TextorLiskiewicz2011,vanconstructing,maathuis2013generalized}
as special cases.

\subsubsection*{Acknowledgements}
This research was supported by the Swiss National Science
Foundation ($200021\_149760$).

\bibliography{biblioteka}

\begin{thebibliography}{}

\bibitem[Ali et~al., 2009]{ali2009markov}
Ali, R.~A., Richardson, T.~S., and Spirtes, P. (2009).
\newblock Markov equivalence for ancestral graphs.
\newblock {\em Ann. Stat.}, 37:2808--2837.

\bibitem[Angrist et~al., 1996]{AngristIR96}
Angrist, J.~D., Imbens, G.~W., and Rubin, D.~B. (1996).
\newblock Identification of causal effects using instrumental variables.
\newblock {\em J. Am. Stat. Assoc.}, 91(434):444--455.

\bibitem[Bareinboim et~al., 2014]{Bareinboim2014}
Bareinboim, E., Tian, J., and Pearl, J. (2014).
\newblock Recovering from selection bias in causal and statistical inference.
\newblock In {\em Proceedings of AAAI 2014}, pages 2410--2416.

\bibitem[Chickering, 2003]{Chickering03}
Chickering, D.~M. (2003).
\newblock Optimal structure identification with greedy search.
\newblock {\em J. Mach. Learn. Res.}, 3:507--554.

\bibitem[Claassen et~al., 2013]{ClaassenMooijHeskes13}
Claassen, T., Mooij, J., and Heskes, T. (2013).
\newblock Learning sparse causal models is not {NP}-hard.
\newblock In {\em Proceedings of UAI 2013}, pages 172--181.

\bibitem[Colombo and Maathuis, 2014]{colombo2014order}
Colombo, D. and Maathuis, M.~H. (2014).
\newblock Order-independent constraint-based causal structure learning.
\newblock {\em J. Mach. Learn. Res.}, 15:3741--3782.

\bibitem[Colombo et~al., 2012]{Colombo2012}
Colombo, D., Maathuis, M.~H., Kalisch, M., and Richardson, T.~S. (2012).
\newblock Learning high-dimensional directed acyclic graphs with latent and
  selection variables.
\newblock {\em Ann. Stat.}, 40:294--321.

\bibitem[Maathuis and Colombo, 2015]{maathuis2013generalized}
Maathuis, M.~H. and Colombo, D. (2015).
\newblock A generalized back-door criterion.
\newblock {\em Ann. Stat.}, 43:1060--1088.

\bibitem[Maathuis et~al., 2010]{MaathuisColomboKalischBuehlmann10}
Maathuis, M.~H., Colombo, D., Kalisch, M., and B\"uhlmann, P. (2010).
\newblock Predicting causal effects in large-scale systems from observational
  data.
\newblock {\em Nat. Methods}, 7:247--248.

\bibitem[Maathuis et~al., 2009]{MaathuisKalischBuehlmann09}
Maathuis, M.~H., Kalisch, M., and B\"uhlmann, P. (2009).
\newblock Estimating high-dimensional intervention effects from observational
  data.
\newblock {\em Ann. Stat.}, 37:3133--3164.

\bibitem[Meek, 1995]{meek1995causal}
Meek, C. (1995).
\newblock Causal inference and causal explanation with background knowledge.
\newblock In {\em Proceedings of UAI 1995}, pages 403--410.

\bibitem[Pearl, 1993]{pearl1993bayesian}
Pearl, J. (1993).
\newblock Comment: Graphical models, causality and intervention.
\newblock {\em Stat. Sci.}, 8:266--269.

\bibitem[Pearl, 2009]{Pearl2009}
Pearl, J. (2009).
\newblock {\em Causality}.
\newblock Cambridge University Press, Cambridge, second edition.

\bibitem[Peters and B\"uhlmann, 2015]{Peters2014neco}
Peters, J. and B\"uhlmann, P. (2015).
\newblock Structural intervention distance ({SID}) for evaluating causal
  graphs.
\newblock {\em Neural Comput.}, 27:771--799.

\bibitem[Richardson and Spirtes, 2002]{richardson2002ancestral}
Richardson, T. and Spirtes, P. (2002).
\newblock Ancestral graph {M}arkov models.
\newblock {\em Ann. Stat.}, 30:962--1030.

\bibitem[Robins, 1986]{robins1986new}
Robins, J. (1986).
\newblock A new approach to causal inference in mortality studies with a
  sustained exposure period-application to control of the healthy worker
  survivor effect.
\newblock {\em Math. Mod.}, 7:1393--1512.

\bibitem[Rubin, 2008]{Rubin2008}
Rubin, D. (2008).
\newblock Author's reply.
\newblock {\em Stat. Med.}, 27:2741--2742.

\bibitem[Shpitser, 2012]{shpitser2012avalidity}
Shpitser, I. (2012).
\newblock Appendum to ``{O}n the validity of covariate adjustment for
  estimating causal effects''.
\newblock Unpublished manuscript.

\bibitem[Shpitser and Pearl, 2006]{shpitser2006identification}
Shpitser, I. and Pearl, J. (2006).
\newblock Identification of joint interventional distributions in recursive
  semi-markovian causal models.
\newblock In {\em Proceedings of AAAI 2006}, pages 1219--1226.

\bibitem[Shpitser et~al., 2010]{shpitser2012validity}
Shpitser, I., VanderWeele, T., and Robins, J.~M. (2010).
\newblock On the validity of covariate adjustment for estimating causal
  effects.
\newblock In {\em Proceedings of UAI 2010}, pages 527--536.

\bibitem[Shrier, 2008]{Shrier2008}
Shrier, I. (2008).
\newblock Letter to the editor.
\newblock {\em Stat. Med.}, 27:2740--2741.

\bibitem[Shrier and Platt, 2008]{ShrierP2008}
Shrier, I. and Platt, R.~W. (2008).
\newblock Reducing bias through directed acyclic graphs.
\newblock {\em BMC Med. Res. Methodol.}, 8(70).

\bibitem[Spirtes et~al., 2000]{spirtes2000causation}
Spirtes, P., Glymour, C., and Scheines, R. (2000).
\newblock {\em Causation, Prediction, and Search}.
\newblock MIT Press, Cambridge, second edition.

\bibitem[Textor and Li\'{s}kiewicz, 2011]{TextorLiskiewicz2011}
Textor, J. and Li\'{s}kiewicz, M. (2011).
\newblock Adjustment criteria in causal diagrams: An algorithmic perspective.
\newblock In {\em Proceedings of UAI 2011}, pages 681--688.

\bibitem[Tian and Pearl, 2002]{tian2002general}
Tian, J. and Pearl, J. (2002).
\newblock A general identification condition for causal effects.
\newblock In {\em Proceedings of AAAI 2002}, pages 567--573.

\bibitem[van~der Zander et~al., 2014]{vanconstructing}
van~der Zander, B., Li\'skiewicz, M., and Textor, J. (2014).
\newblock Constructing separators and adjustment sets in ancestral graphs.
\newblock In {\em Proceedings of UAI 2014}, pages 907--916.

\bibitem[West and Koch, 2014]{West2014}
West, S.~G. and Koch, T. (2014).
\newblock Restoring causal analysis to structural equation modeling.
\newblock {\em Struct. Equ. Modeling}, 21:161--166.

\bibitem[Westreich and Greenland, 2013]{Westreich2013}
Westreich, D. and Greenland, S. (2013).
\newblock {{T}he table 2 fallacy: presenting and interpreting confounder and
  modifier coefficients}.
\newblock {\em Am. J. Epidemiol.}, 177:292--298.

\bibitem[Zhang, 2006]{zhang2006causal}
Zhang, J. (2006).
\newblock {\em Causal Inference and Reasoning in Causally Insufficient
  Systems}.
\newblock PhD thesis, Carnegie Mellon University.

\bibitem[Zhang, 2008]{zhang2008completeness}
Zhang, J. (2008).
\newblock On the completeness of orientation rules for causal discovery in the
  presence of latent confounders and selection bias.
\newblock {\em Artif. Intell.}, 172:1873--1896.

\end{thebibliography}

\end{document}